\documentclass[12pt]{amsart}

\usepackage{amsmath}
\usepackage{hyperref}
\usepackage{xypic}

\newfont{\sheaf}{eusm10 scaled\magstep1}

\def\A{{\mathcal{A}}}
\def\C{{\mathcal{C}}}

\def\T{{\mathfrak{T}}}

\def\hol{{\mathcal{O}}}
\def\M{{\mathcal{M}}}
\def\CC{{\mathbb{C}}}
\def\PP{{\mathbb{P}}}
\def\QQ{{\mathbb{Q}}}
\def\FF{{\mathbb{F}}}
\def\ZZ{{\mathbb{Z}}}
\def\ttt{{\mathfrak{t}}}
\def\PPP{{\mathfrak{P}}}

\def\text{{\hbox}}

\theoremstyle{plain}
\newtheorem{teo}{Theorem}[section]
\newtheorem{lem}[teo]{Lemma}

\newtheorem{prop}[teo]{Proposition}
\newtheorem*{teo*}{Theorem}

\theoremstyle{definition}

\theoremstyle{remark}
\newtheorem{oss}[teo]{Remark}
\newtheorem*{oss*}{Remark}

\DeclareMathOperator{\Aut}{Aut}
\DeclareMathOperator{\Ext}{Ext}
\DeclareMathOperator{\Hom}{Hom}
\DeclareMathOperator{\Pic}{Pic}
\DeclareMathOperator{\Sim}{Sym}
\DeclareMathOperator{\Span}{Span}
\DeclareMathOperator{\rank}{rank}

\begin{document}
\title{Surfaces with $K^2=8$, $p_g=4$ and\\ canonical
  involution.}\thanks{The present work was performed in the realm of
  the DFG-Schwerpunkt {\it Globale Methoden in der komplexen
    Geometrie}. Part of the article was developped while the second
  author was visiting professor in Bayreuth financed by the
  DFG-Forschergruppe {\it Klassifikation algebraischer Fl\"achen und
    kompakter komplexer Mannigfaltigkeiten.}} 

\author[Ingrid C. Bauer]{Ingrid C. Bauer}
\author[Roberto Pignatelli]{Roberto Pignatelli}
\address[Ingrid C. Bauer]{Mathematisches Institut, 
Universit\"atstra\ss e 30, D-95447 Bayreuth, Germany}
\address[Roberto Pignatelli]{Dipartimento di Matematica, Universit\`a
  di Trento, via Sommarive 14 loc. Povo, I-38050 Trento, Italy}

\date{\today}
\maketitle
\markright{SURFACES WITH $K^2=8$, $p_g=4$ AND CANONICAL INVOLUTION}


\section*{Introduction}
The aim of this paper is to classify regular minimal surfaces $S$ with
$K^2 = 8$ and $p_g = 4$ whose canonical map factors through an
involution (short: {\it having a canonical involution}). 

The study of surfaces with geometric genus $p_g = h^0(S, \Omega^2_S)=4$
began with Enriques' celebrated book {\em Le superficie algebriche}
(\cite{enriques}), where he summarized his research of over fifty
years.

By standard inequalities, minimal surfaces with geometric genus $p_g=4$
satisfy $4\leq K_S^2 \leq 45$. While for high values of $K^2_S$ it is
already difficult to prove existence, the challenge for low values is
to completely classify all surfaces with the given value of
$K_S^2$. More ambitiously, one would like to understand the topology
of the moduli space, i.e., the irreducible and connected components of
the moduli space. 

The lowest possible values $K_S^2=4,5$ were already treated by Enriques
and the corresponding moduli spaces were completely understood in the
70's. For $K_S^2=6$ the situation is far more complicated. In
\cite{Hor3} Horikawa completely classifies all surfaces with $p_g=4$
and $K^2=6$, obtaining a stratification of the moduli space in $11$
strata. Moreover he shows that there are $4$ irreducible components,
and at most three connected components. In \cite{Annalen} it is shown
that the number of connected components actually cannot be bigger than
two. Let us point out that all these surfaces are homeomorphic. 

The complete classification of minimal surfaces with $K_S^2=7$ and
$p_g=4$ was achieved by the first author in
\cite{Ingridhabilitation}. Moreover, it is shown there that all these
surfaces are homeomorphic, and that there are three irreducible
components and at most two connected components.

The first open case $K^2_S=8$ is more complicated already for topological
reasons. By work of Ciliberto, Francia, Mendes Lopes, Oliverio and
Pardini (cf. \cite{ciro}, \cite{CFML}, \cite{triple}, \cite{Oliverio})
there are at least three different topological types, therefore at
least three connected components of the moduli space.

The analysis of the cases $K^2 \leq 7$ is based on a detailed study of
the behaviour of the canonical map $\varphi_{K_S}\colon S
\dashrightarrow \PP^3$, as already suggested by Enriques. For $K^2=8$
this approach produces too many strata and the question how they glue
together becomes intractable. Therefore it is necessary to find a less
fine stratification of the moduli space.

We summarize our main result in the following 
\begin{teo*}
Let $S$ be a minimal regular surface with $p_g=4$ and $K^2=8$ whose
canonical map factors through an involution $i$ on $S$. Then: 
\begin{itemize}
\item[1)] the number $\tau$ of isolated fixed points of $i$ is $0,2,4$
  or $20$;
\item[2)] if $\tau=20$, $S$ is a canonical bidouble cover and the two
  additional involutions have $\tau=0$; 
\item[3)] the surface $S$ belongs to exactly one of six unirational
  families. In the table below we give, for each family, the dimension
  and the reference where this family is described; 

\renewcommand{\arraystretch}{1.5}
\begin{tabular}{|l|l|l|l|l|}
\hline
$\mathrm{Family}$ &  $\dim$ & $\mathrm{reference}$ \\
\hline
\hline
$\M_0^{\mathrm{(div)}}$& $29$ & \ref{oliverio}  \\ 
\hline
$\M_0$        & $28$ & \ref{classIV2}  \\ 
\hline
$\M_2^{(0)}$  & $32$ & \ref{classIII20}\\ 
\hline
$\M_2^{(1)}$  & $33$ & \ref{classIII21}\\ 
\hline
$\M_4^{\mathrm{(DV)}}$ & $38$ & \ref{ciliberto} \\
\hline
$\M_4^{(2)}$  & $34$ & \ref{genus2}\\
\hline
\end{tabular}
\item[3)] exactly two of these families, namely $\M_4^{\mathrm{(DV)}}$
  and $\M_4^{(2)}$, are irreducible components of the moduli space;
\item[4)] the surfaces in $\M_0^{\mathrm{(div)}}$ and the surfaces in
  $\M_4^{(2)}$ are not homeomorphic and not homeomorphic to any of the
  others. 
\end{itemize}
\end{teo*}

Actually, we prove more:

\begin{oss*}
The index $\tau \in \{0,2,4\}$ in the above families means that there
is a canonical involution on $S$ having $\tau$ isolated fixed
points. In fact, the only surfaces having more than one canonical
involution are canonical bidouble covers having an involution with
$\tau=20$ and two involutions with $\tau=0$: they give a subfamily of
$\M_0^{\mathrm{(div)}}$ (when the canonical image is a quadric cone)
and a subfamily of $\M_0$ (when the canonical image is smooth).

The surfaces in $\M_0^{\mathrm{(div)}}$ are the only ones in the above
list with $2$-divisible canonical system. The surfaces in
$\M_4^{\mathrm{(DV)}}$ (so called because they are {\it Du Val double
  planes}) are the only ones in this list with nontrivial torsion
subgroup of the Picard group. The surfaces in $\M_4^{(2)}$ are all
minimal surfaces of general type with $K^2=8$ and $p_g=4$ having a
genus $2$ pencil.
\end{oss*}

The paper is organized as follows.

In section $1$ we recall some general facts about involutions and show
that the number of isolated fixed points is $0$, $2$, $4$ or $20$.

Sections $2$, $3$, $4$ and $5$ are devoted to the classification and
to the detailed description of all surfaces having a canonical
involution with respectively $\tau=20$, $\tau=0$, $\tau=2$ and
$\tau=4$.  For $\tau = 0,~2$ we use the MMP for pairs (as e.g. in
\cite{milessmalldegree}). The surfaces (minimal, regular with $p_g=4$
and $K^2=8$) having a canonical involution with $\tau = 4$ are exactly
the surfaces (with the same invariants) whose bicanonical map is not
birational. Those without genus $2$ pencil are classified in
\cite{CFML}. We classify those with a genus $2$ pencil using the
techniques developped in \cite{cp}.

In section $5$ we calculate the dimensions of each family. 

{\bf Acknowledgements:} We are indebted to M. Reid for correctimg a
mistake in a previous version and suggesting to us how to simplify
substantially parts of the paper. 

\section{Canonical involutions}\label{rivdoppi} 

Let $S$ be a  regular minimal surface of general type and let $i$ be
an involution on $S$.

Since $S$ is minimal $i$ is biregular, and its fixed locus consists of
$\tau$ isolated points and a nonsingular (not necessarily connected)
curve $R$.

The quotient $T:=S/i$ has $\tau$ nodes. Resolving them we get a
cartesian diagram of morphisms
\begin{equation}\label{diagrcappuccio}
\xymatrix{
\hat{S}\ar_{\hat{\pi}}[d]\ar^{\epsilon}[r]&S\ar_{\pi}[d]\\
\hat{T}\ar[r]&T
}
\end{equation}
with vertical maps finite of degree $2$ and horizontal maps
birational. We denote by $\Delta$ the branch curve $\pi(R)$ and by
$E_1,\ldots,E_{\tau}$ the exceptional curves of $\epsilon$.

The action of $i$ on $\hat{S}$ yields a decomposition $\hat{\pi}_*
\hol_{\hat{S}}=\hol_{\hat{T}} \oplus \hol_{\hat{T}}(-\hat{\delta})$,
with $2\hat{\delta} \equiv \Delta +\sum_1^\tau \hat{\pi}(E_i)$. Recall
that $K_{\hat{S}}\equiv \hat{\pi}^*(K_{\hat{T}}+\hat{\delta})$.

\begin{lem}\label{chi}
\begin{eqnarray}
\label{chiOT}
\chi(\hol_{\hat{T}})& = & \frac12\chi(\hol_{S})-\frac18(K_SR-\tau),\\
\label{chiOT-delta}
\chi(\hol_{\hat{T}}(-\hat{\delta}))& = &
\frac12\chi(\hol_{S})+\frac18(K_SR-\tau).
\end{eqnarray}
\end{lem}
\begin{proof}
By Riemann-Roch 
\begin{multline*}
\chi(\hol_{\hat{T}})-\chi(\hol_{\hat{T}}(-\hat{\delta}))
=-\frac12\hat{\delta}(K_{\hat{T}}+\hat{\delta})=\\
=-\frac14(R+\sum_1^\tau E_i)(K_S+\sum_1^\tau E_i)=-\frac14(K_SR-\tau).
\end{multline*}
The result follows then from
$\chi(\hol_{S})=\chi(\hol_{\hat{T}})+\chi(\hol_{\hat{T}}(-\hat{\delta}))$.
\end{proof}

We will also use the following (cf. e.g. \cite{rito?})
\begin{equation}\label{rito}
0\leq \tau= K_S^2 + 6\chi(\hol_{\hat{T}}) - 2\chi(\hol_S) -
2h^0(\hol_{\hat{T}}(2K_{\hat{T}} + \hat{\delta})).
\end{equation} 

\begin{oss}
If the canonical map factors through the involution $i$, then either
$p_g(\hat{T})=p_g(S)$ (equivalently, all $2$-forms are invariant) or
$p_g(\hat{T})=0$ (i.e., all $2$-forms are anti invariant).
\end{oss}

\begin{lem}\label{fixedisbase}
Assume that $i$ is a canonical involution and let $p$ be an isolated
fixed point of $i$. 
\begin{itemize}
\item If $p_g(\hat{T})=0$, then $p$ is a base point of $|K_S|$. 
\item If $p_g(\hat{T})=p_g(S)$, then $R$ is contained in the fixed
  part of $|K_S|$.  
\end{itemize}
\end{lem}

\begin{proof}
There are local coordinates around $p$ such that $i(x,y)=(-x,-y)$. 
In particular $i^*(x^ay^b dx \wedge dy)=(-1)^{a+b} x^ay^b dx \wedge dy$. 

If $p_g(\hat{T})=0$, every global $2$-form $\omega$ on $S$ is anti
invariant. Writing $\omega=\sum \omega_{a,b} x^ay^b dx \wedge dy$ it
follows $\omega_{a,b}=0$ for $a+b$ even. In particular $\omega$
vanishes in $p$. 

The other case is similar, since there are local coordinates around
any point of $R$ such that $i(x,y)=(-x,y)$ and $R=\{x=0\}$.
\end{proof}

\begin{oss}\label{derossi}
If $p_g(\hat{T})=0$, by Hurwitz' formula and Riemann-Roch (as in lemma
\ref{chi}) 
\begin{equation*}
\begin{array}{lll}
K_{\hat{T}}\hat{\delta}&=&-2-K^2_{\hat{T}}-\frac12 \tau,\\ 
\hat{\delta}^2&=&8+K^2_{\hat{T}}+\frac12 \tau.
\end{array}
\end{equation*}
\end{oss}

From now on $S$ will be a minimal surface of general type with
$K^2=8$, $p_g=4$, and $q=0$. 

\begin{oss}\label{nopencil}
The canonical map of $S$ is not composed with a pencil.
\end{oss}
More generally, by results of Zucconi and Konno (cf. \cite{zuc} and
\cite{kon}) the canonical map of regular surfaces with  $p_g \geq 3$
and $K^2_S < 4p_g-6$ is not composed with a pencil.

\begin{prop}\label{pg=0}
If the canonical map of $S$ factors through an involution $i$, then either
\begin{itemize}
\item[1)] $p_g(\hat{T})=0$, $\tau \in \{0,2,4\}$, or
\item[2)] $p_g(\hat{T})=4$,  $R=\emptyset$, $\tau = 20$.
\end{itemize}
\end{prop}

\begin{proof}

If $p_g(\hat{T})=4$, the canonical map cannot have degree $2$ (since
then $\hat{T}$ is birational to the canonical image which has degree
at most $4$), therefore it has degree $4$ and $K_S$ is base point
free, so, by lemma \ref{fixedisbase}, $R=\emptyset$. $\tau = 20$
follows from (\ref{chiOT}). 

Otherwise $p_g(\hat{T})=0$. By (\ref{rito})
$\tau=4-2h^0(\hol_{\hat{S}/i}(2K_{\hat{S}/i}+\hat{\delta}))$. 
\end{proof}

\section{Canonical involutions with $p_g(\hat{T})=4$}

In this section $S$ is a minimal surface of general type with $K_S^2=8$, $p_g(S)=4$ and a canonical involution such that $p_g(\hat{T})=4$.

Consider a Hirzebruch surface $\FF_k$ , $k\in \{0,2\}$. Then, if
$k=0$, we denote by $|\Gamma_1|$, $|\Gamma_2|$ the two rulings of
$\FF_0$. Otherwise, we denote by $|\Gamma_2|$ the ruling of $\FF_2$
and by $|\Gamma_1|:=\Gamma_{\infty}+|\Gamma_2|$, $\Gamma_{\infty}$
being the $(-2)$-curve.

We will show the following
\begin{teo}\label{pg=4}
$S$ is a bidouble cover (i.e., a Galois cover with group $\ZZ/2\ZZ
\times \ZZ/2\ZZ$) of $\FF_k$, $k\in \{0,2\}$, which is a fiber product
of two double covers branched in two general divisors
$B_1\in|4\Gamma_1+2\Gamma_2|$, $B_2\in|2\Gamma_1+4\Gamma_2|$.
\end{teo}

First we need the following:
\begin{lem}\label{etale}
Let $C$ be a curve of genus $2$ and let  $f\colon D \rightarrow C$ be
an \'etale double cover with associated involution $\xi$. Then the
hyperelliptic involution of $C$ lifts to an involution on $D$ which
commutes with $\xi$. 
\end{lem}
\begin{proof}
The hyperelliptic involution $\sigma'$ acts on $\Pic^0(C)$ as $L
\mapsto L^*$, and therefore it fixes any $2$-torsion bundle. Since
(connected) \'etale double covers are classified by non trivial
$2$-torsion bundles, considering the fiber product
$$\xymatrix{
D'\ar_{\sigma'^*f}[d]\ar^{\sigma}[r]&D\ar_{f}[d]\\
C\ar^{\sigma'}[r]&C
}
$$
it follows that $D' \cong D$ and $\sigma$ is a lift of $\sigma'$ to $D$.

Since $\sigma'$ is an involution, $\sigma^2$ is either the identity or
$\xi$, which has no fixed points. But in this last case (by Hurwitz)
$D/\sigma$ would have genus $\frac32$, a contradiction.
\end{proof}

By proposition \ref{pg=0}, if $p_g(\hat{T})=4$, then $R=\emptyset$, so
$K^2_T=\frac82=4$. By \cite{Hor1} $T$ is a canonical double cover of
an irreducible quadric in $\PP^3$ branched in the complete
intersection with a general sextic. Moreover the canonical map of $S$
is the composition of $\pi$ with the canonical map of $T$. 

\begin{lem}\label{Galois}
$S$ is a canonical Galois cover of a  quadric in $\PP^3$ with Galois
group $\ZZ/2\ZZ \times \ZZ/2\ZZ$.
\end{lem}
\begin{proof}
The pull-back of a ruling of the quadric is a genus $2$ pencil on $T$
and (since $R=\emptyset$) a genus $3$ pencil on $S$ whose general
element is an \'etale double cover of the corresponding genus $2$
curve. Then by lemma \ref{etale} we can lift the canonical involution
of $T$ to an involution on $S$ commuting with $i$, and the canonical map is
the quotient by these two commuting involutions. 
\end{proof}

$S$ has two more canonical involutions, and we denote them by $\sigma$
and $\sigma i$. 

\begin{lem}\label{sigmatau=0}
$\sigma$ and $\sigma i$ do not have isolated fixed points.
\end{lem}
\begin{proof}
Recall that the action of $i$ on $H^0(K_S)$ is the identity. Since
$p_g(S/_{\ZZ/2\ZZ \times \ZZ/2\ZZ})=0$, the action of $\sigma$ on
$H^0(K_S)$ is multiplication by $-1$, and $p_g(S/\sigma)=p_g(S/\sigma
i)=0$. Since $\deg(\varphi_{|K_S|}) = 4$, $|K_S|$ is base point free,
and the claim follows from lemma \ref{fixedisbase}.
\end{proof}

\begin{proof}[Proof of theorem \ref{pg=4}]
We have a commutative diagram of finite morphisms of degree $2$  
\begin{equation*}
\xymatrix{
&S\ar_{\pi_\sigma}[dl]\ar^{\pi}[d]\ar^{\pi_{\sigma i}}[dr]&\\
S/\sigma\ar_{q_{\sigma i}}[dr]&T\ar^{q}[d]&S/\sigma i\ar^{q_{\sigma}}[dl]\\
&S_{/\ZZ/2\ZZ \times \ZZ/2\ZZ}&
}
\end{equation*}
The ramification locus of $\pi_{\sigma}$ is a smooth divisor
$R_{\sigma}$, the ramification locus of $\pi_{\sigma i}$ is a smooth
divisor $R_{\sigma i}$, the ramification locus of $\pi$ is a set of
$20$ points $\PPP$. $R_{\sigma}$ and $R_{\sigma i}$ intersect
transversally and obviously $\PPP \supset R_{\sigma} \cap R_{\sigma
  i}$. On the other hand, since $\varphi_{|K_T|}$ factors through
$\FF_k$ (\cite{Hor1}, lemma 1.5), the same holds for $\varphi_{|K_S|}$
and therefore $S_{/\ZZ/2\ZZ \times \ZZ/2\ZZ}$ is a blow-up of
$\FF_k$. So  $S_{/\ZZ/2\ZZ \times \ZZ/2\ZZ}$ is smooth, and also the
other inclusion must hold, i.e., $\PPP = R_{\sigma} \cap R_{\sigma i}$. 

We consider the branch divisors  $B_{\sigma}=q \circ \pi (R_{\sigma})$
of $q_{\sigma}$, and $B_{\sigma i}:=q \circ \pi (R_{\sigma i})$ of
$q_{\sigma i}$. It follows that $B_{\sigma} B_{\sigma i}=20$. We
denote by  $D_{\sigma}$, $D_{\sigma i}$ the respective images on
$\FF_k$. Since $B_{\sigma}$, $B_{\sigma i}$ are $2$-divisible, we can
write $B_{\sigma}=D_{\sigma}+2\sum_j a_j E_j$, $B_{\sigma i}=D_{\sigma
  i}+2\sum_j \alpha_j E_j$ where $E_j$ are the exceptional divisors of
the first kind of the contraction to $\FF_k$.

$\pi^*q^* E_i$ is contracted by $\varphi_{|K_S|}$ and $2K_S =
\pi^*q^*( B_{\sigma}+B_{\sigma i}+2K_{S_{/\ZZ/2\ZZ \times
    \ZZ/2\ZZ}})$. Then $a_j+\alpha_j=-1$ for all $j$, so
$a_j\alpha_j\leq 0$ is even and it follows that $D_{\sigma}D_{\sigma
  i}=20-8k$ for some nonnegative integer $k$.

$D_{\sigma}$, $D_{\sigma i}$ are $2$-divisible, effective and
$D_{\sigma}+D_{\sigma i}$ is the branch curve of $q$, so belongs to
$|6\Gamma_1 + 6\Gamma_2|$.  Therefore either $D_{\sigma} \in
|4\Gamma_1 + 2\Gamma_2|$, $D_{\sigma i} \in |2\Gamma_1 + 4\Gamma_2|$,
or $D_{\sigma} \in | 2\Gamma_j|$, $j\in\{1,2\}$. 

A smooth bidouble cover of type (in the language of \cite{bidouble})
($(2,0)$, $(4,6)$, $(0,0)$) has $K^2=8$ and $p_g=6$. By the formulas
on page 109 of \cite{bidouble} there is no configuration of
singularities that changes $p_g$ without changing $K^2$.
\end{proof}

Bidouble covers of a smooth quadric were already studied by Catanese
\cite{bidouble}, and later Gallego and Purnaprajna \cite{gp1} and
\cite{gp2} classified canonical Galois covers of degree $4$ of a
surface of minimal degree. All these surfaces can be found in those
papers. Note however that these surfaces, because of the other two
canonical involutions they have, are also special cases of the
surfaces studied in the next section.

\section{Canonical involutions with $p_g(\hat{T})=0$, $\tau=0$}

In this case $T$ is smooth. By remark \ref{derossi}
\begin{equation}\label{valeanchesuP}
\begin{array}{lll}
K_{T}\delta&=&-2-K^2_{T},\\ 
\delta^2&=&8+K^2_{T}.
\end{array}
\end{equation}

We inductively contract all $(-1)$-curves $E$ on $T$ contained in the
image of the fundamental cycles of $S$, and we denote by $\alpha\colon
T \rightarrow P$ the composition of all these contractions.

\begin{oss}
We observe that every $(-1)$-curve $E$ contained in the image of a
fundamental cycle of $S$ fulfills $\Delta E=2$. It follows that
equations (\ref{valeanchesuP}) hold also for $K_P$, $\delta_P$.
\end{oss}
Let $\lambda \in \mathbb{Q} \cup \{ \infty\}$ be the maximal number
such that $\lambda K_P + \delta_P$ is nef. Since the pull back of $K_P
+ \delta_P$ to $S$ is $K_S$, $\lambda \geq 1$. In fact, $\lambda > 1$,
since $\lambda = 1$ implies that there is an extremal ray $l$ such
that $(K_P + \delta_P)l = 0$. By a.i.t., $l^2 < 0$, whence $l$ is a
$(-1)$-curve whose pull-back to $S$ is contained in a fundamental
cycle. But these have already been contracted.

\begin{prop}\label{KP=01}
There are the following two possibilities:
\begin{itemize}
\item $K_P^2=1$ and $3K_P+\delta_P$ is trivial;
\item $K_P^2=0$ and $|2K_P+\delta_P|$ is a genus $0$ pencil without
  base points.
\end{itemize}
\end{prop}

\begin{proof}
By the algebraic index theorem: $K_P^2\delta_P^2 \leq
(K_P\delta_P)^2$. Equations (\ref{valeanchesuP}) imply $K_P^2 \leq 1$.  

If $K_P^2=1$, equality holds in the a.i.t. and $3K_P+\delta_P$ is
numerically trivial. By equation (\ref{rito}) $2K_P+\delta_P$ is
effective, hence Riemann-Roch implies $h^0(3K_P+\delta_P)\geq
1$. Therefore $3K_P+\delta_P$ is trivial.

Otherwise $K_P^2 \leq 0$. Let $l$ be an extremal ray with $(\lambda
K_P + \delta_P)l = 0$. Since  $P$ is neither $\PP^2$ nor a
$\PP^1$-bundle, $l$ has to be a $(-1)$-curve, whence $\lambda =
\delta_P l \in \ZZ$. In particular, $2K_P + \delta_P$ is nef. 

Since $2K_P+\delta_P$ is effective, whence $ 0 \leq (2K_P +
\delta_P)^2 = K_P^2 \leq 0$. Therefore $2K_P+\delta_P$ is a nef
divisor with selfintersection $0$ and negative canonical degree. This
implies that $|\frac{-2}{K_P(2K_P + \delta_P)} (2K_P + \delta_P)|$ is
a base point free genus $0$ pencil. Since in our case  $K_P(2K_P +
\delta_P)=K_P\delta_P=-2$ we are done.
\end{proof}

We get two families, according to the value of $K_P^2$.

\begin{teo}\label{oliverio}
If $K^2_P=1$, then $K_S$ is $2$-divisible and $S$ is a double cover of
a Del Pezzo surface of degree $1$ branched in a general divisor in
$|{-}6K|$. 
\end{teo}

\begin{proof}
By proposition \ref{KP=01} $3K_P+\delta_P$ is trivial, so
$K_P+\delta_P=-2K_p$ is $2$-divisible and the same holds for its
pull-back $K_S=\pi^*\alpha^*(K_P+\delta_P)$. Note that since
$K_P+\delta_P$ is ample, $P$ is a Del Pezzo surface.
\end{proof}

\begin{oss}\label{oliveriodescription}
Oliverio proves in \cite{Oliverio} that if the canonical system of a
regular minimal surface with $K^2_S=8$ and $p_g=4$ is $2$-divisible,
either $K_S$ has base points and the canonical map has degree $3$ (so
it is not our case), or the semicanonical ring $R(S,\frac12 K_S)$
embeds the canonical model of $S$ as a complete intersection of two
sextics in $\PP(1,1,2,3,3)$. 
\end{oss}

\begin{teo}\label{classIV2}
If $K_P^2=0$, then $S$ is the minimal resolution of a double cover of
$\PP^1 \times \PP^1$ branched in a curve in $|8\Gamma_1+10\Gamma_2|$
having $8$ singular (possibly infinitely near) of multiplicity $4$ as
only essential singularities.
\end{teo}

\begin{proof}
By proposition \ref{KP=01}, $|2K_P+\delta_P|$ is a genus $0$ pencil
without base points. 

Contracting $8-K_{P}^2=8$ ($-1$)-curves (contained in fibres) we get a
birational morphism $\eta\colon P\rightarrow \FF_r$. Note that here
one has to choose the $8$ contractions and different choices yield
different $r$'s. We show that we can choose the contractions so that $r=0$.

Else, the strict transform of the $(-r)$-section $\Gamma_{\infty}$ of
$\FF_r$ is an irreducible rational curve $B_{\infty}$ on $P$ with
$B_{\infty}(2K_P+\delta_P)=1$. Let $E$ be a $(-1)$-curve contained in
a fibre of $|2K_P+\delta_P|$: then $B_{\infty}E$ is $0$ or $1$. If
$B_{\infty}E=1$, $2K_{P}+\delta_P-E$ is again an exceptional divisor
of the first kind, so it contains an other $(-1)$-curve $E'$ and
$B_{\infty}E'=0$. Therefore we can choose $\eta$ such that for all
contracted curves holds $B_{\infty}E=0$.

Now, $B_{\infty}$ is a smooth rational curve with $B_{\infty}^2=-r$,
so $K_PB_{\infty}=r-2$. Therefore $0 \leq (K_P+\delta_P)B_{\infty}=
(2K_P+\delta_P)B_{\infty}-K_PB_{\infty}=3-r$ whence $r\leq 3$.

We modify our choice of $\eta$. We contract exactly $r$ $(-1)$-curves
with $B_{\infty}E=1$ in order to get $\eta\colon P\rightarrow \PP^1
\times \PP^1$. 

We choose the rulings of $\PP^1 \times \PP^1$ such that
$\eta^*\Gamma_2=2K_P+\delta_P$. We write
$\delta_P=\eta^*(a\Gamma_1+b\Gamma_2)-\sum_1^8c_iE_i$. 

First of all, for all $i$,
$c_i=\delta_PE_i=(2K_P+\delta_P)E_i-2K_PE_i=2$.  Moreover, by formulae
(\ref{valeanchesuP}) $a= \delta_P(\eta^*\Gamma_2)= \delta_P(2K_P+\delta_P)= -4+8= 4$. Finally,  by $8=\delta_P^2=8b-32$ we get $b=5$, so $\delta_P=\eta^*(4\Gamma_1+5\Gamma_2)-2\sum_1^8E_i$.
\end{proof}

\section{Canonical involutions with $p_g(\hat{T})=0$, $\tau=2$}

We recall diagram (\ref{diagrcappuccio})
\begin{equation*}
\xymatrix{
\hat{S}\ar_{\hat{\pi}}[d]\ar^{\epsilon}[r]&S\ar_{\pi}[d]\\
\hat{T}\ar[r]&T
}
\end{equation*}

In this case $\epsilon$ is the blow up of $S$ in two distinct points
$p_1$ and $p_2$. We denote by $A_i$ the $(-2)$-curve
$\hat{\pi}(\epsilon^{-1}(p_i))$. Note that $A_i$ is a component of the
branch curve of $\hat{\pi}$ with $\hat{\delta}A_i=-1$.

We define the $\mathbb{Q}$-divisor $\bar{\delta} := \hat{\delta} -
\frac{1}{2} (A_1 + A_2)$. We have 
\begin{equation}\label{valeanchesuT1}
\begin{array}{lll}
K_{\hat{T}}\bar{\delta}&=&-3-K^2_{\hat{T}},\\ 
\bar{\delta}^2&=&10+K^2_{\hat{T}}.
\end{array}
\end{equation}

Observe that, by a.i.t., $K_{\hat{T}}^2 (10 + K_{\hat{T}}^2) =
K_{\hat{T}}^2 \bar{\delta}^2 \leq (K_{\hat{T}}\bar{\delta})^2 =
(K_{\hat{T}}^2 + 3)^2$, therefore $K_{\hat{T}}^2 \leq 2$.

Let $\lambda$ be the maximal (rational) number such that $\lambda
K_{\hat{T}} + \bar{\delta}$ is nef. Note that $\hat{\pi}^*(K_{\hat{T}}
+ \bar{\delta}) = \epsilon^* K_S$, whence $K_{\hat{T}} + \bar{\delta}$
is nef, so $\lambda \geq 1$.

Assume that $\lambda = 1$ and let $l$ be an extremal ray with
$(K_{\hat{T}} + \bar{\delta})l = 0$. Since $K_{\hat{T}}^2 \leq 2$, we
know that $\hat{T}$ is neither $\PP^2$ nor a $\PP^1$-bundle. Therefore
$l$ is a $(-1)$-curve, and we contract it. Note that after this
contraction the equations (\ref{valeanchesuT1}) remain valid (if by
slight abuse of notation we denote the pushforward of $\bar{\delta}$
again by $\bar{\delta}$), since $K^2$, $\bar{\delta}^2$ increase by
$1$, while $K \bar{\delta}$ decreases by $1$. In particular, by the
index theorem we get $K^2 \leq 2$.

Therefore, we can inductively apply the above argument and get a
sequence of contractions $c_1 \colon \hat{T} \rightarrow P$, such that
(\ref{valeanchesuT1}) holds on $P$ (so $K_P^2 \leq 2$) and there are
no extremal rays in $(K_{P}+ \bar{\delta})^{\perp}$.

Now, let $\lambda$ be the maximal rational number such that $\lambda
K_{P}+ \bar{\delta}$ is nef. Then $\lambda > 1$.

Since $K_{P}^2 \leq 2$, an extremal ray $l$ has to be a $(-1)$-curve,
whence $\lambda = \bar{\delta} l \in \frac{1}{2} \ZZ$ (since
$2\bar{\delta}$ is integral),  i.e., $\lambda \geq \frac{3}{2}$. 

In particular, $\frac{3}{2} K_{P}+ \bar{\delta}$ is nef and, since by
(\ref{rito}) $2K_{P} + \bar{\delta}$ is effective, we have $0 \leq
(\frac{3}{2} K_{P}+ \bar{\delta})(2K_{P} + \bar{\delta}) =
\frac{1}{2}(K_{P}^2 - 1)$. Therefore, $K_{P}^2 \in \{1, 2\}$.

\begin{prop}\label{pignarompe}
One of the following occurs:
\begin{itemize}
\item $K_{P}^2=2$ and $|4K_{P}+2\bar{\delta}|$ is a genus $0$ pencil
  without base points; 
\item $K_{P}^2=1$, there is a birational morphism $c \colon P
  \rightarrow P_1$ onto a Del Pezzo surface of degree $5$, contracting
  $(-1)$-curves $l$ with $(K + \bar{\delta})l = \frac{1}{2}$, and
  $2K_{P_1} + \bar{\delta} \equiv 0$. 
\end{itemize}
\end{prop}

\begin{proof}
We know that $\lambda \geq \frac{3}{2}$. Assume that $\lambda =
\frac{3}{2}$ and let $l$ be an extremal ray with $(\frac{3}{2}K +
\bar{\delta})l = 0$. By a.i.t., since $(\frac{3}{2}K + \bar{\delta})^2
= 1 + \frac{K^2}{4} > 0$, $l^2 < 0$. Contracting $l$ we add $1$ to
$K^2$, $\frac{9}{4}$ to $\bar{\delta}^2$, and we subtract
$\frac{3}{2}$ from $K \bar{\delta}$, in particular, we do not change
$(\frac{3}{2}K + \bar{\delta})^2$. Therefore we can repeat the
argument and inductively contract all, say $s$, $(-1)$-curves $l$ with
$(\frac{3}{2}K + \bar{\delta})l = 0$. We get a birational morphism $c
\colon P \rightarrow P_1$, such that on $P_1$, $\lambda > \frac32$. 

Since $K_{P_1} + \bar{\delta}$ is nef and $2K_{P_1} + \bar{\delta}$ is
effective, we have $0 \leq (K_{P_1} + \bar{\delta})(2K_{P_1} +
\bar{\delta}) = 1 - \frac{s}{4}$, i.e., $s \leq 4$. In particular,
$K_{P_1}^2 \leq 2 + s \leq 6$, so, as above, an extremal ray $l$ has
to be a $(-1)$-curve, whence $\lambda = \bar{\delta} l \in \frac{1}{2}
\ZZ$, so $\lambda \geq 2$. Therefore $2K_{P_1} + \bar{\delta}$ is a
nef and effective divisor with selfintersection $(2K_{P_1} +
\bar{\delta})^2=(2K_{P} +
\bar{\delta})^2+\frac{s}{4}=K_P^2-2+\frac{s}{4}$.

If $K_{P}^2=1$, it follows that $s=4$ and $(2K_{P_1} +
\bar{\delta})^2=(K_{P_1} + \bar{\delta})(2K_{P_1} +
\bar{\delta})=0$. By a.i.t. $2K_{P_1} + \bar{\delta}$ is trivial.

Else $K_{P}^2=2$, and the inequality $K_{P_1}^2 \bar{\delta}^2 \leq
(K_{P_1}\bar{\delta})^2$ gives $(2+s)(12+\frac94s) \leq (5+\frac32
s)^2 \Leftrightarrow s \leq \frac23$: we have $s=0$. In this case,
$P=P_1$ and $2K_P+\delta_P$ is nef with selfintersection $0$ and
canonical degree $-1$, so $|2(2K_P+\bar{\delta}))|$ is a base point
free genus $0$ pencil. 
\end{proof}

Therefore we get two families, according to the value of $K_P^2$.

\begin{teo}\label{classIII20}
If $K_P^2 = 1$, then  $S$ is the minimal resolution of a double cover
of a Del Pezzo surface of degree $5$ branched in a divisor in
$|{-}4K|$ having two $(3,3)$ points.
\end{teo}

\begin{proof}
Let $l \subset \hat{T}$ be a $(-1)$-curve with $(K_{\hat{T}} +
\bar{\delta})l = 0$. Since the intersection form restricted to
$(K_{\hat{T}}+ \bar{\delta})^{\perp}$ is negative definite and since
$l$, $A_1$, $A_2 \in (K_{\hat{T}}+ \bar{\delta})^{\perp}$, $l(A_1+A_2)
\leq 1$. Because $\bar{\delta}l$, $\hat{\delta}l \in  \mathbb{Z}$, we
have $l(A_1+A_2)$ even, thus $lA_1 = lA_2 =0$. 

This shows that the images of $A_1$ and $A_2$ are still $(-2)$-curves
in $P$. We show that they will be contracted by $c$. Recall that $c$
is (any) sequence of $4$ contractions of extremal rays in $(\frac32
K_{\hat{T}}+ \bar{\delta})^{\perp}$. 

The first extremal ray is a $(-1)$-curve $l$ with $(\frac32
K_{\hat{T}} + \bar{\delta})l = 0$. By the same argument as above,
$l(A_1+A_2) \leq 1$. $\bar{\delta}l \not\in \mathbb{Z}$,
$\hat{\delta}l \in  \mathbb{Z}$, therefore w.l.o.g. $lA_1=1$ and $lA_2
=0$.

After contracting $l$, $A_1$ becomes a $(-1)$-curve contained in
$(\frac32 K_{\hat{T}}+ \bar{\delta})^{\perp}$, and we can choose $A_1$
as second extremal ray. 

By the same argument the third extremal ray $l'$ has $l'A_2 = 1$ and
we can choose $A_2$ as last extremal ray.

Now, $P_1$ is a Del Pezzo of degree $5$, and $P$ is the blow up of
$P_1$ in four points. We call the exceptional divisors $E_1, \ldots,
E_4$.  By the above arguments, we can assume that $A_1 = E_3 - E_4$,
$A_2 = E_1 - E_2$, and on $P$ we have: $\bar{\delta} =c^*\bar{\delta}
- \sum_{i=1}^4 (\bar{\delta}E_i)E_i = c^*(-2K_{P_1}) - \sum_{i=1}^4
\frac32 E_i$.  The direct image of $\hat{\delta}$ on $P$ is therefore
$c^*(-2K_{P_1}) - \sum_{i=1}^4 2E_i + E_1 + E_3$.
\end{proof}

\begin{teo}\label{classIII21}
If $K_P^2=2$, then $S$ is the minimal resolution of a double cover of
$\PP^1 \times \PP^1$ branched in a fibre $\Gamma \in |\Gamma_2|$ and a
curve in $|8\Gamma_1+9\Gamma_2|$ having $6$ singular points $x_1,
\ldots, x_6$ of multiplicity $4$ as only essential singularities, with
$x_5 \in \Gamma$ and $x_6$ infinitely near to $x_5$ and belonging to
the strict transform of $\Gamma$. 
\end{teo}

\begin{proof}
By proposition \ref{pignarompe}, $|4K_P+2\bar{\delta}|$ is a genus $0$
pencil without base points. As in the previous proof we note that
$A_1$ and $A_2$ are still $(-2)$-curves on $P$, which are contained in
fibres of the pencil. 

Contracting $8-K_{P}^2=6$ $(-1)$-curves (contained in fibres) we get a
birational morphism $\eta\colon P\rightarrow \FF_r$. Repeating the
same argument as in the proof of \ref{classIV2}, we obtain that $r
\leq \frac52$ and that therefore we can choose the $6$ exceptional
curves appropriately to get $\eta\colon P\rightarrow \PP^1 \times \PP^1$. 

Let $l$ be one of these $6$ $(-1)$-curves; then being $l, A_1$ and
$A_2$ all contained in fibres, by Zariski's lemma $lA_i \leq 1$. But
it cannot be $lA_1=0$ for all $l$, since $\PP^1 \times \PP^1$ does not
contain curves with negative selfintersection. Therefore one of these
extremal rays has $lA_1=1$, say $E_6$. $\bar{\delta}l$, $\hat{\delta}l
\in  \mathbb{Z}$, therefore $E_6(A_1+A_2)$ even, thus $E_6A_1 = E_6A_2
=1$. So, after this contraction $A_1$ and $A_2$ become ($-1$)-curves
contained in a fibre with $A_1A_2=1$. One will be contracted and the
other will map isomorphically onto a line of $\PP^1 \times \PP^1$. 

We choose the rulings of $\PP^1 \times \PP^1$ such that
$\eta^*\Gamma_2=4K_P+2\bar{\delta}$. We write
$\bar{\delta}=\eta^*(a\Gamma_1+b\Gamma_2)-\sum_1^6c_iE_i$. 

Then, for all $i$, $c_i=\bar{\delta}E_i=2$.  Moreover, by formulae
(\ref{valeanchesuT1})
$a=\bar{\delta}(\eta^*\Gamma_2)=\bar{\delta}(4K_P+2\bar{\delta})=4$.
Finally, by $12=\bar{\delta}^2=8b-24$ we get $b=\frac92$, so
$\bar{\delta}=\eta^*(4\Gamma_1+\frac92 \Gamma_2)-2\sum_1^6E_i$.

Therefore the direct image of $\hat{\delta}$ on $P$ is
$\eta^*(4\Gamma_1+5 \Gamma_2)-\sum_1^5 2E_i - 3E_6$.
\end{proof}

\begin{oss}
We observe that the surfaces classified in this section are exactly
those whose canonical map is a double cover of a cubic surface in $\PP^3$. 
\end{oss}

\section{Canonical involutions with $p_g(\hat{T})=0$, $\tau=4$}

This case can be treated with the same techniques as in the previous
two sections, but the calculations become more demanding. We choose a
different approach. 

By equation (\ref{rito}),
$h^0(\hol_{\hat{S}/i}(2K_{\hat{S}/i}+\hat{\delta}))=0$, in particular,
the bicanonical map factors through the involution $i$. In \cite{CFML}
the authors classify all surfaces with $p_g \geq 4$, nonbirational
bicanonical map having no genus $2$ pencil. In particular, they obtain 

\begin{teo}[\cite{CFML}, thm. 3.1 and rem. 3.10]\label{ciliberto} 
If $\tau=4$ and $S$ has no genus $2$ pencil, then $S$ belongs to one
of the following two families 

\begin{itemize}
\item[i)] $S$ is birational to a double cover of $\PP^1 \times \PP^1$
  with branch curve $\tilde{\Delta}=L_1+L_1'+L_2+L_2'+D$ where
  $L_i,L_i'$ are distinct lines in $|\Gamma_i|$ and $D\in |8 \Gamma_1
  + 8 \Gamma_2|$ has quadruple points at the intersection of the $4$
  lines as only essential singularities.
\item[ii)] $S$ is birational to a double cover of $\PP^1 \times \PP^1$
  with branch curve $\tilde{\Delta}=L_2+L_2'+D$ where $L_2,L_2'$ are
  distinct lines in $|\Gamma_2|$ and $D\in |8 \Gamma_1 + 8 \Gamma_2|$
  has $(4,4)$ points at the intersection of the $2$ lines with a line
  $L_1$ in $|\Gamma_1|$, having as tangent line $L_2$ resp. $L_2'$, as
  only essential singularities.
\end{itemize}
The torsion subgroup of $\Pic(S)$ is isomorphic to $\ZZ/2\ZZ$. The
second case is a specialization of the first one.
\end{teo}

\begin{oss}
It is well known that, if a surface has a genus $2$ pencil, the
involution on each fibre induces an involution on $S$ such that both
the canonical and the bicanonical map of $S$ factor through it. In
particular, the induced involution is canonical and, if the surfaces
is regular with $K^2=8$ and $p_g=4$, it has $\tau=4$.  

It follows that none of the preceedingly studied surfaces has a genus $2$
pencil. 
\end{oss}

In the following $S$ is assumed to be a surface of general type with
$K^2=8$, $p_g=4$ and $q=0$ having a genus $2$ pencil $f\colon S
\rightarrow \PP^1$.

\begin{oss}
Since $\tau = 4$, the canonical system has base points (cf. lemma
\ref{fixedisbase}) and therefore the canonical map has degree two onto
a cubic or a quadric.
\end{oss}

Let $\omega_{S|\PP^1}:=\omega_S \otimes f^*\omega_{\PP^1}^{-1}$ be the
relative canonical sheaf. The sheaves $f_*\omega_{S|\PP^1}^n$ are
vector bundles and there are the relative n-canonical maps
$\varphi_n\colon S \dashrightarrow
\PP(f_*\omega_{S|\PP^1}^n):=\mathbf{Proj}(\Sim
f_*\omega_{S|\PP^1}^n)$, whose restriction to each fibre is its
n-canonical map. Note that for $g=2$ the target of the relative
n-canonical map is a $\PP^1$-bundle for $n=1$ and a $\PP^2$-bundle for
$n=2$.

\begin{oss}\label{linetoline}
Let $f\colon S \rightarrow \PP^1$ be a genus $2$ fibration with fibres
$f^{-1}(t)=:F_t \in |F|$ and assume
\begin{equation}\label{noncompen}
\forall t \in \PP^1 \hbox{ the restriction map } H^0(\omega_S)
\rightarrow H^0(\omega_{F_t}) \hbox{ is surjective.}
\end{equation}
Then the canonical map of $S$ factors through the relative canonical
map. The resulting map $\PP(f_*\omega_{S|\PP^1}) \rightarrow
\varphi_{|K_S|}(S)$ is a surjective morphism mapping each ``line'' of
the ruling of $\PP(f_*\omega_{S|\PP^1})$ to a line of $\PP^{p_g-1}$.

If $S$ is regular, then the cokernels of the restriction maps in
(\ref{noncompen}) are all isomorphic (to $H^1(\omega_S(-F))$). In
particular, the maps are {\em all} surjective if and only if one of
them is surjective, i.e., if and only if $|K_S|$ is not composed with $|F|$.
\end{oss}

\begin{oss}
The canonical map of $S$ is a double cover of a quadric. In fact, by
the above considerations the canonical image is covered by lines. On
the other hand, as it is seen by the same argument as in lemma 3.14 of
\cite{Ingridhabilitation}, if the canonical image of $S$ is a cubic,
it has isolated singularities, whence cannot be covered by lines.
\end{oss}

\begin{prop}\label{twocases}
Let $S$ be a regular surface, whose canonical map is a double cover of
a quadric surface $Q$, and let $f\colon S \rightarrow \PP^1$ be a
genus $2$ fibration. If $Q$ is smooth then $f_*\omega_{S|\PP^1}\cong
2\hol_{\PP^1}(3)$. If $Q$ is a quadric cone then
$f_*\omega_{S|\PP^1}\cong \hol_{\PP^1}(2) \oplus \hol_{\PP^1}(4)$.
\end{prop} 

\begin{proof} 
$\PP(f_*\omega_{S|\PP^1})$ is a Hirzebruch surface $\FF_k$ having, by
remark \ref{linetoline}, a birational morphism onto $Q$. If $Q$ is
smooth, then $k=0$, and if the quadric is a cone, then $k=2$. We
conclude, since by standard computations (e.g., \cite{cp}, rem. 2.11)
$\deg f_*\omega_{S|\PP^1}=\chi(\hol_S)+1=6$.
\end{proof}

\begin{lem}\label{nobranchinfty}
With the same hypotheses as in proposition \ref{twocases}, if $Q$ is a
quadric cone, then the branch curve of the relative canonical map
$\varphi_1\colon S \dashrightarrow \FF_2$ cannot contain $\Gamma_{\infty}$.
\end{lem}
\begin{proof}
Assume by contradiction that $\Gamma_{\infty}$ is contained in the
branch locus of $\varphi_1$. Then the preimage of the vertex of the
cone under the canonical map is a point $p \in S$. Since the genus two
pencil maps onto the ruling of $Q$, it has a base point, contradicting
Kodaira's lemma (\cite{HorP} or \cite{xiao}, prop. 5.1).
\end{proof}

We will use some of the techniques developped in \cite{cp}, which for
sake of simplicity will only be briefly reported in the case of genus
$2$ fibrations $f\colon S \rightarrow \PP^1$ with $p_g(S)=4$. 

We consider the exact sequence 
\begin{equation}\label{sigma2}
0
\rightarrow
\Sim^2 f_*\omega_{S|\PP^1}
\stackrel{\sigma_2}{\rightarrow}
f_*\omega_{S|\PP^1}^2
\rightarrow
\hol_{\ttt}
\rightarrow
0,
\end{equation}
where $\sigma_2$ is the natural map induced by the tensor product of
canonical sections of the fibers of $f$, and $\ttt$ is an effective
divisor on $\PP^1$ of degree $K^2_S-4$ (cf. lemma 4.1 of
\cite{cp}). The map $\sigma_2$ yields a rational map
$\nu\colon\PP(f_*\omega_{S|\PP^1})\dashrightarrow
\PP(f_*\omega_{S|\PP^1}^2)$ (relative version of $2$-Veronese
embedding $\PP^1 \hookrightarrow \PP^2$) birational onto a conic
bundle $\C$. 

The following exact sequence defines the vector bundle $\A_6$ as
quotient of $\Sim^3 f_*\omega_{S|\PP^1}^2$, the vector bundle of
relative cubics on $\PP(f_*\omega_{S|\PP^1}^2)$, by the subbundle of
cubics vanishing on $\C$ (cf. lemma 4.4 of \cite{cp}):
\begin{equation}\label{A6}
0
\rightarrow
f_*\omega_{S|\PP^1}^2 \otimes \hol_{\PP^1}(12)
\stackrel{i_3}{\rightarrow}
\Sim^3 f_*\omega_{S|\PP^1}^2
\rightarrow
\A_6
\rightarrow
0.
\end{equation}
The branch curve $\Delta$ of the map $S \rightarrow \C$ is given
(cf. thm. 4.7 and prop. 4.8 of  \cite{cp}) by a map
\begin{equation}\label{Delta}
\delta \colon \hol_{\PP^1}(2K_S^2+4) \hookrightarrow \A_6.
\end{equation}

\begin{lem}\label{atleast6}
Under the assumptions of proposition \ref{twocases}, if moreover
$K^2_S \geq 6$, then each direct summand of $f_*\omega_{S|\PP^1}^2$
has degree at least $6$. 
\end{lem}
\begin{proof}
Being $\sigma_2$ an injective morphism between two vector bundles of
the same rank, if each summand of the source has degree at least $6$,
the same holds for the target. Therefore by prop. \ref{twocases} we
can assume $f_*\omega_{S|\PP^1}=\hol_{\PP^1}(2)\oplus
\hol_{\PP^1}(4)$.

Assume by contradiction that (writing coordinates on
$f_*\omega_{S|\PP^1}$, $f_*\omega_{S|\PP^1}^2$) 
\begin{equation*}
 f_*\omega_{S|\PP^1}=x_0\hol_{\PP^1}(2)\oplus x_1\hol_{\PP^1}(4)
\end{equation*}
\begin{equation*}
 f_*\omega_{S|\PP^1}^2=y_0\hol_{\PP^1}(a) \oplus y_1\hol_{\PP^1}(b)
 \oplus y_2\hol_{\PP^1}(c)
\end{equation*}
with $a\leq 5$. 

In these coordinates we have that $\Gamma_{\infty}$ has equation
$x_1=0$. From $a\leq 5$ it follows that $\sigma_2(x_0x_1)$,
$\sigma_2(x_1^2)$ belong to $\Span(y_1,y_2)$, whence
$\nu(\Gamma_{\infty})=\{y_1=y_2=0\}$.

Since $\nu(\Gamma_{\infty}) \subset \C$,  $y_0^2$ does not appear in
the equation of $\C$ and therefore $y_0^3$ does not appear in the
equation of any relative cubic vanishing in $\C$. This means that the
row of the matrix of $i_3$ corresponding to the direct summand
$y_0^3\hol_{\PP^1}(3a)$ of $\Sim^3 f_*\omega_{S|\PP^1}^2$ is a line of
zeroes. Therefore this summand maps isomorphically onto a direct
summand of $\A_6$.

$K_S^2 \geq 6$ implies $2K_S^2+4>15\geq 3a$ and therefore the
composition of $\delta$ with the projection on this summand is
zero. But this implies $\Delta \supset \nu(\Gamma_{\infty})$,
contradicting lemma \ref{nobranchinfty}.
\end{proof}

Let now $S$ be a minimal surface of general type with $K^2=8$, $p_g=4$
and $q=0$ having a genus $2$ pencil $f\colon S \rightarrow \PP^1$.

By the above arguments we know:

\begin{itemize}
\item $\PP(f_*\omega_{S|\PP^1})\cong \FF_k$ for $k\in\{0,2\}$;
\item $f_*\omega_{S|\PP^1}^2\cong r\hol_{\PP^1}(6) \oplus V$ for 
$r\in \{0,1,2\}$, where $V$ is a sum of line
bundles of degree at least $7$.
\end{itemize}

Note that $r \neq 3$, since $\deg f_*\omega_{S|\PP^1}^2=18+\deg \ttt
=22$.

\begin{teo}\label{genus2}
The moduli space of surfaces with $K^2=8$, $p_g=4$ and $q=0$ having a
genus $2$ pencil $f\colon S \rightarrow \PP^1$ is unirational of
dimension $34$.
\end{teo} 

\begin{proof}
We use the structure theorem for genus $2$ fibrations (cf. thm.
4.13 in \cite{cp}). For each case we have to describe the associated 
$5$-tuple $(B,V_1,\ttt,\xi,w)$. We treat separately the cases $k=0$
and $k=2$. 

\medskip

$\mathbf{k=0.}$ The first three elements are easy: $B=\PP^1$,
$V_1=f_*\omega_{S|\PP^1}=2\hol_{\PP^1}(3)$ and $\ttt$ is an effective
divisor on $\PP^1$ of degree $4$.

$\xi$ is an element of $\Ext^1_{\hol_{\PP^1}}(\hol_{\ttt},
\Sim^2V_1)_{/\Aut_{\hol_{\PP^1}}(\hol_\ttt)}$, giving the short exact
sequence (\ref{sigma2}). In order to give explicitly these extension
classes we fix a section $f_{\ttt}\in H^0(\hol_{\PP^1}(\ttt))$ and,
applying to the exact sequence 
\begin{equation}\label{37}
0
\rightarrow
\hol_{\PP^1}(3)
\stackrel{f_{\ttt}}{\rightarrow}
\hol_{\PP^1}(7)
\rightarrow
\hol_{\ttt}
\rightarrow
0
\end{equation}
the functor $\Hom_{\hol_{\PP^1}}(\cdot,3\hol_{\PP^1}(6))$, we get 
$$\Ext^1_{\hol_{\PP^1}}(\hol_{\ttt},
\Sim^2V_1) \cong \Hom_{\hol_{\PP^1}}(\hol_{\PP^1}(3),3\hol_{\PP^1}(6)) \cong
H^0(3\hol_{\PP^1}(3))\cong \CC^{12}.$$

This isomorphism is explicitly given as follows: for any triple of
cubics $(c_0,c_1,c_2)$, the resulting $f_*\omega_{S|\PP^1}^2$ is given
by the short exact sequence
\begin{equation}
0
\rightarrow\label{37666}
\hol_{\PP^1}(3)
\stackrel{c}{\longrightarrow}
\hol_{\PP^1}(7)\oplus 3\hol_{\PP^1}(6)
\rightarrow
f_*\omega_{S|\PP^1}^2
\rightarrow
0
\end{equation}
for $c$ being the transpose of $(-f_{\ttt},c_1,c_2,c_3)$; $\sigma_2$
is then the restriction to the last three summands
($3\hol_{\PP^1}(6)$) of the projection on $f_*\omega_{S|\PP^1}^2$.

These $4$ data give us the exact sequence (\ref{sigma2}) and therefore
the conic bundle $\C$. To complete the $5$-tuple we have to give an
element $w\in (\Hom(\hol_{\PP^1}(20),\A_6)\setminus \{0\})_{/\CC^*}$
corresponding to the map $\delta$ in (\ref{Delta}), and then to the
branch curve $\Delta \subset \C$.

From the exact sequence (\ref{A6}), $\dim(\Hom(\hol_{\PP^1}(20),\A_6)) 
=\chi(\A_6(-20))+h^1(\A_6(-20))=29+h^1(\A_6(-20))$. Moreover,
$H^1(\A_6(-20))$ is isomorphic to the cokernel of the map $H^1(i_3(-20))$. 

By lemma \ref{atleast6}, all summands of the source and of the target
of the map $i_3(-20)$ have degree at least $-2$. More precisely, the
source has $r\leq 2$ summands of degree $-2$, the target $r^2$, and
$H^1(i_3(-20))$ is a map $\CC^r \rightarrow \CC^{r^2}$. In particular, 
\begin{equation}\label{randh1}
r^2-r\leq h^1(\A_6(-20)) \leq r^2.
\end{equation}
In fact, the map $H^1(i_3(-20))$ is easily obtained by the matrix of
$i_3$ by taking the $r^2 \times r$ submatrix $A$ given by the rows and
the columns of the summands of degree $18$ (both in the source and in
the target).

We have three cases, according to the value of $h^1(\A_6(-20))$. 

{$\underline{h^1(\A_6(-20))=0}$}: This happens for a general choice of
$\xi$, since dualizing the exact sequence (\ref{37666}) one sees that,
if the three cubics $c_1,c_2,c_3$ are linearly independent, $r=0$.

We have $4$ parameters for $\ttt$, $12-4=8$ for $\xi$ and $29-1=28$
for $w$: 40 parameters. Since we must take the quotient by the action
of $\Aut(\PP^1 \times \PP^1)$, this family is unirational of dimension
$34$.

{$\underline{h^1(\A_6(-20))=1}$}: By (\ref{randh1}), then
$r=1$, i.e., there is a nontrivial relation $\alpha c_1+\beta
c_2+\gamma c_3=0$ between the three cubics: these are two conditions
for $\xi$. Moreover, the row of the matrix of $\sigma_2$ corresponding
to the degree $6$ summand of the target is $(\alpha, \beta, \gamma)$,
and $A=(\alpha\gamma-\beta^2)$. In order to get $r=1$ we need to
further assume $\alpha\gamma=\beta^2$; we have three conditions on
$\xi$, and therefore this gives a family of dimension
$34-3+h^1(\A_6(-20))=32$. 

{$\underline{h^1(\A_6(-20))\geq 2}$}: By (\ref{randh1}), then $r=2$,
i.e., the three cubics span a space of dimension $1$: these are six
conditions.  Moreover, if the submatrix of $\sigma_2$ corresponding to
the degree $6$ summands of the target is
\begin{equation}\label{2x3}
\begin{pmatrix}
\alpha_1&\beta_1&\gamma_1\\
\alpha_2&\beta_2&\gamma_2
\end{pmatrix},
\end{equation}
the matrix $A$ is
$$\begin{pmatrix}
\alpha_1\gamma_1-\beta_1^2&0\\
\alpha_1\gamma_2+\alpha_2\gamma_1-2\beta_1\beta_2&\alpha_1\gamma_1-\beta_1^2\\
\alpha_2\gamma_2-\beta_2^2&\alpha_1\gamma_2+\alpha_2\gamma_1-2\beta_1\beta_2\\
0&\alpha_2\gamma_2-\beta_2^2
\end{pmatrix}.$$
It follows: $\rank A \neq 2 \Leftrightarrow A =0$. If $A=0$, then
$(\alpha_1y_1 + \alpha_2y_2)(\gamma_1y_1 + \gamma_2y_2)-(\beta_1y_1 +
\beta_2y_2)^2 = 0$, and this implies that the matrix (\ref{2x3}) has
not rank $2$, contradicting the injectivity of $\sigma_2$. Therefore
$h^1(\A_6(-20))=2$ and this gives a family of dimension
$34-6+h^1(\A_6(-20))=30$. 

\medskip

$\mathbf{k=2}.$ Here $V_1=\hol_{\PP^1}(2)\oplus \hol_{\PP^1}(4)$. The
main difference to the first case is that here to describe the
extension class we need to apply the functor
$\Hom_{\hol_{\PP^1}}(\cdot,\Sim^2V_1)$ to the exact sequence 
\begin{equation}\label{15}
0
\rightarrow
\hol_{\PP^1}(1)
\stackrel{f_{\ttt}}{\rightarrow}
\hol_{\PP^1}(5)
\rightarrow
\hol_{\ttt}
\rightarrow
0
\end{equation}
getting only a short exact sequence
\begin{multline*}
0
\rightarrow
\Hom(\hol_{\PP^1}(5),\Sim^2V_1)
\rightarrow\\
\rightarrow
\Hom(\hol_{\PP^1}(1),\Sim^2V_1)
\rightarrow
\Ext^1(\hol_\ttt,\Sim^2V_1)
\rightarrow
0.
\end{multline*}
To induce any extension as described in the first case we need maps
$\hol_{\PP^1}(1) \rightarrow \hol_{\PP^1}(4) \oplus \hol_{\PP^1}(6)
\oplus\hol_{\PP^1}(8)$ (but not in a unique way): the dimension of the
$\Ext^1$ is in fact $18-6=12$ as in the first case. We distinguish two
cases.

{$\underline{h^1(\A_6(-20))=0}$}: This happens for general choice of
$\xi$, since also in this case, if $\xi$ is general, then $r=0$. The
analysis of this case is identical to the analogous case for $k=0$, so
we find again $40$ parameters. Since $\dim \Aut(\FF_2) = 7$, we get an
unirational family of dimension $40-7=33$. 

{$\underline{h^1(\A_6(-20))\geq 1}$}: By (\ref{randh1}) in this case 
$r\geq 1$. Let us first assume $r=1$: then the row of the matrix of 
$\sigma_2$ corresponding to the degree $6$ summand of the target is
$(\alpha, \beta, 0)$ (where $\deg \alpha=2, \beta\in \CC$), and
therefore the matrix $A$ is $(-\beta^2)$. It follows that
$h^1(\A_6(-20))= 1$ forces $\beta=0$.

We are now in the same situation as in the proof of lemma
\ref{atleast6}: $\sigma_2(x_0x_1)$, $\sigma_2(x_1^2)$ belong to $\Span
(y_1,y_2)$. Arguing as there, we conclude that $\Delta \supset
\nu(\Gamma_{\infty})$ contradicting lemma \ref{nobranchinfty}.

The case $r=2$ is similar and even easier, since in this case we can
always assume (up to a change of coordinates in the target) that the
submatrix of $\sigma_2$ corresponding to the degree $6$ summands has
the form 
$$
\begin{pmatrix}
\alpha_1&0&0\\
\alpha_2&\beta_2&0
\end{pmatrix}.
$$

Summing up we have found 4 families, one generically smooth
unirational of dimension $34$, say the ``main'' family, and three more
of respective dimensions $32$, $30$ and $33$. To conclude, we have to
show that the general surface in each of those last three families
admits a small deformation to a surface belonging to the "main"
family.

This is easy for surfaces in the family with $k=2$. In fact, we first
deform $\FF_2$ to $\FF_0$ (i.e., the vector bundle $V_1$). Then,
leaving $\ttt$ fixed, we can deform the extension class $\xi$, since
all the $\Ext^1$ groups have the same dimension $12$: geometrically
this corresponds to deform $\C$ to a family of conic bundles. Finally,
we can deform the last datum, $w$, since we have seen that (for $k=2$)
$h^1(\A_6(-20))=0$, so by semicontinuity it must be zero also on
nearby fibres, and therefore $h^0(\A_6(-20))$ remains constant for a
small deformation: this geometrically corresponds to deform $\Delta$.

This argument does not work for the other two families, since in these
cases $h^1(\A_6(-20))\neq 0$ and therefore, once we have fixed a
$1$-parameter deformation of $\C$, we will not be able to deform all
possible curves $\Delta$.

We use a different argument. Each of the two families is contained in
a irreducible component of the subscheme of the moduli space given by
the surfaces having a canonical involution. We claim that it has
dimension at least $34$.

For the general surface in each of our two families, $\C$ has $\deg
\ttt=4$ nodes (the vertices of the singular conics), none of them in
$\Delta$, which is smooth. Let $\tilde{\C}$ be a minimal
desingularization of $\C$; the $4$ $(-2)$-curves on $\tilde{\C}$ give
rise to $4$ $(-1)$-curves on the associated double cover $\tilde{S}$,
the exceptional locus of the birational morphism $\tilde{S}
\rightarrow S$. The finite double cover $\varphi\colon \tilde{S}
\rightarrow \tilde{\C}$ branches in $\tilde{\Delta}$, union of the
pull-back of $\Delta$ with the $(-2)$-curves.

The invariant part of $\varphi_* (\Omega^1_{\tilde{S}} \otimes
\Omega^2_{\tilde{S}})$ is isomorphic to $\Omega^1_{\tilde{\C}}(\log
\tilde{\Delta}) \otimes \Omega^2_{\tilde{\C}}$.  

The morphism $\tilde{\C} \rightarrow \PP(V_1)$ is the contraction of
the strict transforms of each component of the singular conics, so of
$2\deg \ttt=8$ exceptional curves of the first kind. If $\T_{\cdot}$
denotes the tangent sheaf,
$\chi(\T_{\tilde{C}})=\chi(\T_{\FF_0})-4\deg \ttt=6-16=-10$. Then our
claim follows from 
\begin{multline*}
h^1(\Omega^1_{\tilde{\C}}(\log \tilde{\Delta}) \otimes
\Omega^2_{\tilde{\C}})- 
h^2(\Omega^1_{\tilde{\C}}(\log \tilde{\Delta}) \otimes
\Omega^2_{\tilde{\C}})\geq\\
\geq 
-\chi(\Omega^1_{\tilde{\C}}(\log \tilde{\Delta}) \otimes
\Omega^2_{\tilde{\C}})=
-\chi(\Omega^1_{\tilde{\C}}\otimes \Omega^2_{\tilde{\C}})
-\chi(\hol_{\tilde{\Delta}}(\Omega^2_{\tilde{\C}})))=\\
=-\chi(\T_{\tilde{\C}})-\chi(\Omega^2_{\tilde{\C}})+\chi(\Omega^2_{\tilde{\C}}(-\tilde{\Delta})) 
=10+\frac12 \tilde{\Delta}(\tilde{\Delta}-K_{\tilde{\C}})=\\
=6+\frac12 \Delta(\Delta-K_{\C})=34
\end{multline*}
where $\Delta(\Delta-K_{\C})=56$ is a standard intersection
computation (note that $\C \in |\hol_{\PP(V_2)}(2) \otimes
\hol_{\PP^1}(-12)|$, $\Delta$ is a divisor in the linear system
induced on $\C$ by $|\hol_{\PP(V_2)}(3) \otimes \hol_{\PP^1}(-20)|)$.

Then, since for a small deformation preserving the involution also the
bicanonical map factors through it, either the two families are in the
closure of the "main" family or these surface can be deformed to
surfaces  as in thm. \ref{ciliberto}. But this is impossible for
topological reasons, since the surfaces in thm. \ref{ciliberto} have
non trivial $2$-torsion in $\Pic(S)$ whereas every surface with a
linear pencil of genus $2$ curves and slope $<3$ (in our case
$\frac83$) is simply connected by \cite{Xiaopi1}, theorem 3.
\end{proof}

\section{Moduli}

In the previous sections we classified all pairs $(S,i)$ where $S$ is a
minimal regular surfaces with $K_S^2 =8$, $p_g=4$, and $i$ is a canonical
involution on $S$, finding $8$ families.

\renewcommand{\arraystretch}{1.3}
\begin{tabular}{|l|c|l|}
\hline
{\it Family}& {\it Theorem}& {\it short description}\\
\hline \hline
$\M_+^{(0)}$& \ref{pg=4}&  bidouble covers of $\FF_0$ branched in
  $(4,2)$, $(2,4)$\\
\hline
$\M_+^{(2)}$& \ref{pg=4}& bidouble covers of $\FF_2$ branched in
  $(4,2)$, $(2,4)$\\
\hline
$\M_0^{\mathrm{(div)}}$
& \ref{oliverio}& double covers of a Del Pezzo of degree $1$\\
&& branched in $-6K$\\
\hline
$\M_0$ & \ref{classIV2} & double covers of $\FF_0$\\
&& branched in $(8,10)-4\sum_1^8 E_i$\\
\hline
$\M_2^{(0)}$& \ref{classIII20} &double covers of a Del Pezzo of degree $5$\\
&& branched in $-4K$ with two $(3,3)$\\
\hline 
$\M_2^{(1)}$& \ref{classIII21} & double covers of $\FF_0$\\
&& branched in $(8,10)$ with certain singularities\\ 
\hline
$\M_4^{(2)}$&\ref{genus2}& the surfaces having a genus $2$ pencil\\
\hline
$\M_4^{\mathrm{(DV)}}$ &\ref{ciliberto}& $2K$ non birational, but no
  genus $2$ pencil\\
\hline
\end{tabular}

\begin{oss}
The first two are the families for which $H^0(K_S)$ is
invariant. These surfaces have in fact (lemma \ref{Galois} and
\ref{sigmatau=0}) two more involutions for which $H^0(K_S)$ is
anti invariant and $\tau=0$. In fact, for the family $\M_+^{(0)}$ the
two further involutions are in $\M_0$, for the family $\M_+^{(2)}$ the
two further involutions are in $\M_0^{\mathrm{(div)}}$.

On the other hand, since the canonical map has maximal degree $4$, if
one of these surface has more than one canonical involution, it must
have one involution for which $H^0(K_S)$ is invariant: so these two
families give all surfaces having more than one canonical involution.
\end{oss}

Our results yield then a stratification of the corresponding subscheme
of the moduli space of minimal regular surfaces of general type with
$K_S^2 =8$, $p_g=4$ in six families, image of the last $6$ families of
the above table.

\medskip

The aim of this section is to prove the following

\begin{teo}\label{unirationality}
$\M_4^{\mathrm{(DV)}}$ and $\M_4^{(2)}$ give unirational irreducible
components of the moduli space of minimal regular surfaces of general
type with $K_S^2 =8$, $p_g=4$ of respective dimensions $38$ and $34$. 

The remaining $4$ families $\M_0^{\mathrm{(div)}}$, $\M_0$, $\M_2^{(0)}$,
$\M_2^{(1)}$ give unirational strata of respective dimensions
$29,28,32,33$.
\end{teo}

\begin{oss}
By Kuranishi's theorem each irreducible component of the moduli space
of minimal surfaces of general type with $K^2=8$, $p_g=4$ has
dimension at least $10\chi-2K^2=34$. It follows that the last four
families are not irreducible components of the moduli space.

Observe that the general point of the irreducible component in which
each of these families is contained is a surface without a canonical
involution. In fact, it cannot be in $\M_4^{\mathrm{(DV)}}$ or in
$\M_4^{(2)}$ because $\tau$ is invariant under deformations preserving
the involution.
\end{oss}

\begin{oss}
$\M_4^{\mathrm{(DV)}}$ and $\M_4^{(2)}$ are generically smooth. This
is proved in \cite{Supino} for $\M_4^{\mathrm{(DV)}}$. The same
calculation as in \cite{Ingridhabilitation}, theorem 5.32, shows it
for $\M_4^{(2)}$. 
\end{oss}

\begin{oss}
Minimal surfaces of general type with $K^2=8$, $p_g=4$ belong to at
least three different topological types (in particular, the moduli
space has at least three connected components). The surfaces in
$\M_0^{\mathrm{(div)}}$ are the only ones in our list with
$2$-divisible canonical class, the surfaces in  $\M_4^{\mathrm{(DV)}}$
are the only ones in our list with non trivial torsion in the Picard
group.
\end{oss}

\begin{proof}[Proof of theorem \ref{unirationality}] 
The statement about $\M_4^{\mathrm{(DV)}}$ is theorem 1.3 of \cite{Supino}. 

By theorem \ref{genus2} $\M_4^{(2)}$ is unirational of dimension
$34$. To prove that it is an irreducible component of the moduli space
we need to show that for a general surface in this family the
anti invariant part (with respect to the involution) of
$H^1(\Omega^1_{S} \otimes \Omega^2_{S})$ is trivial. 

This computation works almost identically as the analogous one in
\cite{Ingridhabilitation}, section 5.3. We sketch it.

Using the same notation as in the proof of theorem \ref{genus2},
recall that for a general surface $S$ in $\M_4^{(2)}$, we have a
finite double cover $S \rightarrow \C=S/i$ branched in the $\deg
\ttt=4$ nodes of $\C$, and in the smooth divisor $\Delta$. Resolving
the singular points of $\C$ and blowing up their preimages in $S$ we
get a finite double cover $\varphi\colon \tilde{S} \rightarrow
\tilde{\C}$ whose branch locus is a smooth divisor $\tilde{\Delta}$,
union of the pull-back of $\Delta$ with the $(-2)$-curves. 

Now we can compute the  dimension of the anti invariant part of
$H^1(\Omega^1_{\tilde{S}}\otimes \Omega^2_{\tilde{S}})$ with respect
to the lifting of the involution $i$ to $\tilde{S}$ exactly as in the
proof  of theorem 5.32  of \cite{Ingridhabilitation}: the result is
$8$. Since $b\colon \tilde{S} \rightarrow S$ is a sequence of $4$ blow
ups, by lemma 5.34 of \cite{Ingridhabilitation} the dimension of the
anti invariant part of $H^1(\Omega^1_{S}\otimes \Omega^2_{S})$ is
$8-2\cdot 4=0$.

\smallskip

We prove now the second part of the statement. In all $4$ cases $S$ is
a double cover of a surface $P$ such that the movable part of the
branch curve is $2\delta$ where $\delta$ is a $\QQ$-divisor such that
$\lambda K_P +\delta$ is ample for $\lambda \leq 1$. In particular,
$2\delta_P-K_P$ is ample, therefore $h^1(2\delta_P)=0$, and the
dimension of the linear system $|2\delta|$ can be computed by
Riemann-Roch.

\smallskip

$\underline{\M_0^{\mathrm{(div)}}}$: Del Pezzo surfaces of degree $1$
are obtained by choosing $8$ points in $\PP^2$, therefore, modulo
$\Aut(\PP^2)$, they depend on $8$ (unirational) parameters. Curves in
$|{-}6K|$ depend on $1+\frac12 (42K^2)-1=21$ parameters.

$\underline{\M_0}$: $P$ is the blow up of $\FF_0$ in $8$ general
points, branched in a curve in
$|8\Gamma_1+10\Gamma_2-4\sum_1^8E_i|$. Since $8$ points in $\FF_0$
depend on $16$ parameters and $\dim \Aut(\FF_0)=6$, $P$ depends on
$10$ parameters. The branch curve depends on $18$ parameters.

$\underline{\M_2^{(0)}}$: $P$ is the blow up of a Del Pezzo of degree
$5$ in $4$ points, $2$ of which are infinitely near to the other
two. Therefore $P$ depends on $6$ parameters. The branch curve depends
on $26$ parameters.

$\underline{\M_2^{(1)}}$: $P$ is the blow up of $\FF_0$ in $6$ points,
the last determined by the previous one. Therefore $P$ depends on
$10-6=4$ parameters. The branch curve depends on $29$ parameters.
\end{proof}


\begin{thebibliography}{10}


\bibitem[Bau]{Ingridhabilitation} Bauer, I.C.
{\it Surfaces with $K^2=7$ and $p_g=4$.}
 Mem. Amer. Math. Soc.  152  (2001),  no. 721.

\bibitem[BCP]{Annalen} 
Bauer, I.C.; Catanese, F.; Pignatelli, R.
{\it The moduli space of surfaces with $p_g=4$ and $K^2=6$.}
Math. Ann.  336  (2006),  no. 2, 421--438.




\bibitem[Cat]{bidouble} 
Catanese, F.
{\it Singular bidouble covers and the construction of interesting algebraic surfaces.}     Algebraic geometry: Hirzebruch 70 (Warsaw, 1998),  97--120, 
Contemp. Math., 241, Amer. Math. Soc., Providence, RI, 1999.

\bibitem[Cil]{ciro}
Ciliberto, C.
{\it Canonical surfaces with $p_g=p_a=4$ and $K_S^2=5,\ldots,10$.}
Duke Math. J.  48  (1981), no. 1, 121--157.

\bibitem[CFML]{CFML} Ciliberto, C.; Francia, P.; Mendes Lopes, M.
{\it Remarks on the bicanonical map for surfaces of general type.}
 Math. Z.  224  (1997),  no. 1, 137--166.

\bibitem[CP]{cp} Catanese, F.; Pignatelli R.
{\it Fibrations of low genus, I.}
Ann. Sci. \'Ecole Norm. Sup. (4)  39  (2006),  no. 6, 1011--1049.

\bibitem[Enr]{enriques}
Enriques, F.
{\it Le superficie algebriche}.
Nicola Zanichelli, Bologna, 1949.

\bibitem[GP1]{gp1}
Gallego F.J.,  Purnaprajna B.P.
{\it Classification of quadruple Galois canonical covers. I.}
preprint math/0302045.

\bibitem[GP2]{gp2}
Gallego F.J.,  Purnaprajna B.P.
{\it Classification of quadruple Galois canonical covers. II.}
J. Algebra 312 (2007), no. 2, 798--828.


\bibitem[Hor1]{Hor1} 
Horikawa, E.
{\it Algebraic surfaces of general type with small $c^2_1$. I.}
 Ann. of Math. (2)  104  (1976), no. 2, 357--387.

\bibitem[Hor2]{HorP} 
Horikawa, E.
{\it On algebraic surfaces with pencils of curves of genus $2$.}
Complex analysis and algebraic geometry,  pp. 79--90.
Iwanami Shoten, Tokyo, 1977.

\bibitem[Hor3]{Hor3} 
Horikawa, E.
{\it Algebraic surfaces of general type with small $c^2_1$. III.}
Invent. Math.  47  (1978), no. 3, 209--248.

\bibitem[Kon]{kon} 
Konno, K.
{\it On the irregularity of special non canonical surfaces.} 
Publ. Res. Inst. Math. Sci.  30  (1994),  no. 4, 671--688.

\bibitem[MLP1]{triple}
Mendes Lopes, M., Pardini, R.
{\it Triple canonical surfaces of minimal degree.}
Internat. J. Math.  11  (2000),  no. 4, 553--578.

\bibitem[MLP2]{rito?}
Mendes Lopes M., Pardini R. ´
{\it The bicanonical map of surfaces with $p_g=0$ and $K^2\geq 7$. II.}
Bull. London Math. Soc.  35  (2003),  no. 3, 337--343.

\bibitem[Oli]{Oliverio}
Oliverio, P.A.
{\it On even surfaces of general type with $K^2=8,~p_g=4,~q=0$.}
Rend. Sem. Mat. Univ. Padova  113  (2005), 1--14.

\bibitem[Rei]{milessmalldegree} 
Reid, M.
{\it Surfaces of small degree.}
Math. Ann. 275 (1986), no. 1, 71--80.

\bibitem[Sup]{Supino} 
Supino, P.
{\it On moduli of regular surfaces with $K^2=8$ and $p_g=4$},
Port. Math. (N.S.) 60, no. 3, 353--358 (2003).

\bibitem[Xia1]{xiao} 
Xiao, G.
{\it Surfaces fibr\'ees en courbes de genre deux.}
Lecture Notes in Mathematics, 1137. Springer-Verlag, Berlin, 1985.

\bibitem[Xia2]{Xiaopi1} 
Xiao, G. 
{\it $\pi\sb 1$ of elliptic and hyperelliptic surfaces.}
Internat. J. Math.  2  (1991),  no. 5, 599--615.

\bibitem[Zuc]{zuc}
Zucconi, F.
{\it Numerical inequalities for surfaces with canonical map composed with a
pencil.} 
Indag. Math. (N.S.)  9  (1998),  no. 3, 459--476.

\end{thebibliography}
\end{document}